\documentclass[11pt]{article}
\usepackage{amsmath,amssymb,amsthm,color}
\usepackage[bookmarks,colorlinks]{hyperref}
\usepackage{natbib}

\newtheorem{theorem}{Theorem}

\theoremstyle{definition}

\newcommand{\set}[1]{\left\{#1\right\}}
\newcommand{\paren}[1]{\left(#1\right)}

\newcommand{\bracq}[1]{\left[#1\right]_q}

\begin{document}



\centerline{\Large{\bf Quadratic addition rules for three $q$-integers }}

\centerline{}

\centerline{Mongkhon Tuntapthai}

\centerline{Department of Mathematics, Faculty of Science, Khon Kaen University}

\centerline{Khon Kaen 40002, Thailand}

\centerline{}

\centerline{mongkhon@kku.ac.th}

\begin{abstract}
The $q$-integer is the polynomial $\bracq n = 1 + q + q^2 + \dots + q^{n-1}$. 
For every sequences of polynomials 
$\mathcal S = \set{s_m(q)}_{m=1}^\infty$, 
$\mathcal T = \set{t_m(q)}_{m=1}^\infty$,
$\mathcal U = \set{u_m(q)}_{m=1}^\infty$ and
$\mathcal V = \set{v_m(q)}_{m=1}^\infty$, 
define an addition rule for three $q$-integers by 
\[
	\oplus_{\mathcal S,\mathcal T,\mathcal U,\mathcal V}
	(\bracq m, \bracq n, \bracq k)
	= s_m (q) \bracq m + t_m (q) \bracq n + u_m(q) \bracq k + v_m (q) \bracq n \bracq k .
\]
This is called the first kind of quadratic addition rule for three $q$-integers, if 
\[\oplus_{\mathcal S,\mathcal T,\mathcal U,\mathcal V} 
(\bracq m, \bracq n, \bracq k) 
= \bracq{m+n+k}\]
for all positive integers $m$, $n$, $k$. 

In this paper the first kind of quadratic addition rules for three $q$-integers are determined when $s_m(q)\equiv 1$. Moreover, the solution of the functional equation for a sequence of polynomials $\set{f_n(q)}_{n=1}^\infty$ given by 
\[
	f_{m+n+k} (q) = f_m (q) + q^m f_n (q) + q^m f_k (q) + q^m (q-1) f_n (q) f_k (q)
\] 
for all positive integers $m$, $n$, $k$, are computed.
\end{abstract}





\section{Introduction}

For every positive integer $n$, the $q$-integer $\bracq n$ is the polynomial
\[
	\bracq n = 1 + q + q^2 + \dots +q^{n-1}
\]
and define $\bracq 0=0$. These polynomials appear in many contexts such as in quantum calculus \cite{kac2001quantum}, and quantum groups \cite{kassel2012quantum}.

From the fact that $\bracq{m+n} = \bracq m + q^m \bracq n$ for all positive integers $m$ and $n$, \cite{nathanson2007linear} 
defined the linear addition rule for $q$-integers by
\begin{align*}
	\bracq m \oplus_\ell \bracq n = \bracq m + q^m \bracq n
\end{align*}
for all positive integers $m$ and $n$. Associated to the linear addition rule for $q$-integers, the functional equation for a sequence of polynomials $\set{f_n(q)}_{n=1}^\infty$, given by
\begin{align*}
	f_{m+n} (q) = f_m(q) + q^m f_n(q) ,
\end{align*}
were also studied. The solution of this euqation are defined by $f_n(q) = \bracq n f_1(q)$. 

From two observations that $\bracq{m+n} = \bracq m + \bracq n + (q-1) \bracq m \bracq n$ and that $\bracq{m+n} = \bracq m + \bracq n + (q-1) \bracq m \bracq n$ for all positive integers $m$ and $n$, \cite{kontorovich2006quadratic} defined two non-linear addition rules for $q$-integers, respectively, by 
\begin{align*}
	\bracq m \oplus_1 \bracq n 
	& = \bracq m + \bracq n + (q-1) \bracq m \bracq n ,
	\\
	\bracq m \oplus_2 \bracq n 
	& = q^n \bracq m + q^m \bracq n + (1-q) \bracq m \bracq n ,
\end{align*}
for all positive integers $m$ and $n$. 
These give the corresponding two functional equations on a sequence of polynomials $\set{f_n(q)}_{n=1}^\infty$, given by
\begin{align*}
	f_{m+n} (q) 
	& = f_m(q) + f_n(q) + (q-1) f_m(q) f_n(q) ,
	\\
	f_{m+n} (q) 
	& = q^n f_m(q) + q^m f_n(q) + (1-q) f_m(q) f_n (q) ,
\end{align*}
whose solutions are defined, respectively, by
\begin{align*}
	f_n (q) 
	& = \frac{1-\set{1+(q-1)f_1(q)}^n}{1-q} ,
	\\
	f_n (q) 
	& = \frac{\set{q+(1-q)f_1(q)}^n-q^n}{1-q} .
\end{align*}

Moreover, \cite{nathanson2003functional} defined the multiplication rule for $q$-integers and considered the associated polynomial functional equation. Many articles contained solutions of the multiplicative functional equations, for examples \cite{nguyen2013classification, nguyen2012polynomial, nguyen2010solutions, nguyen2010support, nathanson2004formal}.

In this paper, we consider an addition rule for three $q$-integers which is a mixed type of linear addition rule and the first kind of quadratic addition rule, given by
\begin{align*}
	\oplus \paren{\bracq m , \bracq n, \bracq k}
	= \bracq m \oplus_\ell \paren{ \bracq n \oplus_1 \bracq k } 
\end{align*}
for all positive integers $m$, $n$, $k$, or equivalently,
\begin{align*}
	\oplus \paren{\bracq m , \bracq n, \bracq k}
	= \bracq m + q^m \bracq n + q^m \bracq k + q^m (q-1) \bracq n \bracq k .
\end{align*}
Furthermore, we shall compute the solution of the functional equation 
\begin{align*}
	f_{m+n+k} (q) 
	& = f_m(q) + q^m f_n(q) + q^m f_k(q) + q^m (q-1) f_n(q) f_k(q)
\end{align*}
on a sequence of polynomials $\set{f_n(q)}_{n=1}^\infty$. 



\section{Quadratic addition rules}

For every sequences of polynomials 
$\mathcal S = \set{s_m(q)}_{m=1}^\infty$, 
$\mathcal T = \set{t_m(q)}_{m=1}^\infty$,
$\mathcal U = \set{u_m(q)}_{m=1}^\infty$ and
$\mathcal V = \set{v_m(q)}_{m=1}^\infty$, 
define an addition rule for three $q$-integers by 
\begin{align*}
	\oplus_{\mathcal S,\mathcal T,\mathcal U,\mathcal V} \paren{\bracq m, \bracq n, \bracq k} = s_m (q) \bracq m + t_m (q) \bracq n + u_m(q) \bracq k + v_m (q) \bracq n \bracq k .
\end{align*}
This is called the first kind of quadratic addition rule for three $q$-integers, if 
\begin{align*}
	\oplus_{\mathcal S,\mathcal T,\mathcal U,\mathcal V} \paren{\bracq m, \bracq n, \bracq k} = \bracq{m+n+k}
\end{align*}
for all positive integers $m$, $n$ and $k$. 

\begin{theorem}
\label{thm:quaruleam}
For every sequence of constants 
$\mathcal A = \set{a_m}_{m=1}^\infty$, 
and three sequences of polynomials 
$\mathcal T = \set{t_m (q)}_{m=1}^\infty$, 
$\mathcal U = \set{u_m (q)}_{m=1}^\infty$, 
$\mathcal V = \set{v_m (q)}_{m=1}^\infty$, 
the $\oplus_{\mathcal A,\mathcal T,\mathcal U,\mathcal V}$ is the first kind of quadratic addition rule for three $q$-integers if and only if $a_m=1$, $t_m (q) = u_m (q) = q^m$, and $v_m (q)=q^m(q-1)$ for all positive integers $m$.
\end{theorem}

\proof
Let $\oplus_{\mathcal A,\mathcal T,\mathcal U,\mathcal V}$ be the first kind of quadratic addition rule for three $q$-integers. Thus 
\begin{align}
	\bracq{m+n+k} 
	= a_m \bracq m + t_m (q) \bracq n + u_m (q) \bracq k + v_m (q) \bracq n \bracq k
	\label{eq:coef:1}
\end{align}
for all positive integers $m,n,k$. 
For all positive integers $m$ and $n$, we have
\begin{align*}
	\bracq{m+n+1} = a_m \bracq m + t_m (q) \bracq n + u_m (q) \bracq 1 + v_m (q) \bracq n \bracq 1
\end{align*}
and
\begin{align*}
	\bracq{m+1+n} = a_m \bracq m + t_m (q) \bracq 1 + u_m (q) \bracq n + v_m (q) \bracq 1 \bracq n
\end{align*}
Subtracting, we obtain
\begin{align*}
	0 = t_m (q) \paren{ \bracq n - \bracq 1} - u_m (q) \paren{\bracq n - \bracq 1} ,
\end{align*}
which imples that $t_m (q) = u_m (q)$ for all $m\ge1$. 
For all positive integers $m$, we have 
\begin{align*}
	\bracq{m+1+1} = a_m \bracq m + t_m (q) \bracq 1 + u_m (q) \bracq 1 + v_m (q) \bracq 1 \bracq 1
\end{align*}
and
\begin{align*}
	\bracq{m+1+2} = a_m \bracq m + t_m (q) \bracq 1 + u_m (q) \bracq 2 + v_m (q) \bracq 1 \bracq 2
\end{align*}
Subtracting, we obtain
\begin{align*}
	q^{m+2} = q u_m (q) + q v_m (q) ,
\end{align*}
it follows that $v_m (q) = q^{m+1} - u_m (q)$ for all $m$. 
Now, the equation (\ref{eq:coef:1}) can be expanded to
\begin{align*}
	\bracq{m+n+k} = a_m \bracq m + u_m (q) \bracq n + u_m (q) \bracq k + \paren{q^{m+1}-u_m(q)} \bracq n \bracq k
\end{align*}
for all positive integers $m,n,k$. 
For all positive integers $m$, we obtain 
\begin{align*}
	\bracq{m+1+1} = a_m \bracq m + u_m (q) \bracq 1 + u_m(q) \bracq 1 + \paren{q^{m+1} -u_m (q)} \bracq 1 \bracq 1 .
\end{align*}
From this fact, we can show that 
\begin{align*}
	\bracq{m+1+1} = a_m \bracq m + u_m (q) + q^{m+1} 
\end{align*}
and so that 
\begin{align*}
	\bracq{m+1} = a_m \bracq m + u_m (q) .
\end{align*}
Then, $u_m (q) = q^m + (1-a_m) \bracq m$ and $v_m (q) = q^m (q-1) + (a_m-1) \bracq m$.
By replacing $u_m(q)$ and $v_m(q)$ in (\ref{eq:coef:1}), 
\begin{align*}
	\bracq{m+n+k} 
	& = a_m \bracq m 
	+ q^m \bracq n + (1-a_m) \bracq m \bracq n 
	+ q^m \bracq k + (1-a_m) \bracq m \bracq k 
	\\
	& \quad
	+ q^m(q-1) \bracq n \bracq k + (a_m-1) \bracq m \bracq n \bracq k 
	\\
	& = 
	\bracq m + q^m \bracq n + q^m \bracq k 
	+ q^m(q-1) \bracq n \bracq k
	\\
	& \quad
	+ (a_m-1) \bracq m \set{ 1 - \bracq n - \bracq k + \bracq n \bracq k}
	\\
	& = 
	\bracq{m+n+k} + (a_m-1) \bracq m \set{ 1 - \bracq n - \bracq k + \bracq n \bracq k} ,
\end{align*}
we can see that $0 = (a_m-1) \bracq m \set{ 1 - \bracq n - \bracq k + \bracq n \bracq k}$ for all $m,n,k\ge1$.
Therefore $a_m=1$, and hence $t_m(q)=u_m(q)=q^m$,  $v_m(q)=q^m(q-1)$ for all positive integers $m$.
\endproof


\begin{theorem}
For every sequences of polynomials 
$\mathcal S = \set{s_m}_{m=1}^\infty$, 
$\mathcal T = \set{t_m (q)}_{m=1}^\infty$, 
$\mathcal U = \set{u_m (q)}_{m=1}^\infty$, 
$\mathcal V = \set{v_m (q)}_{m=1}^\infty$, 
if there exists a positive integer $m$ such that the degree of a polynomial $s_m(q)$ is greater than $2$, then $\oplus_{\mathcal S,\mathcal T,\mathcal U,\mathcal V}$ is not the first kind of quadratic addition rule for three $q$-integers.
\end{theorem}

\proof
Suppose that $\deg s_m(q)>2$ for some positive integer $m\ge1$, and that $\oplus_{\mathcal S,\mathcal T,\mathcal U,\mathcal V}$ is the first kind of quadratic addition rule for three $q$-integers. Then
\begin{align*}
	\bracq{m+n+k} = s_m (q) \bracq m + t_m (q) \bracq n + u_m (q) \bracq k + v_m (q) \bracq n \bracq k
\end{align*}
for all positive integers $m,n,k$.
In the proof of Theorem \ref{thm:quaruleam}, it is easy to see that 
$t_m (q)  = u_m (q)$, $v_m (q) = q^{m+1} - u_m (q)$, and $u_m(q) = \bracq{m+1} - s_m (q) \bracq m$ for all $m\ge1$. So, we get
\begin{align*}
	\deg u_m(q)
	& = \deg s_m(q) \bracq m > m+1 ,
	\\
	\deg u_m(q)
	& > \deg s_m(q) ,
	\\
	\deg v_m(q) 
	& = \deg u_m(q) >m+1 ,
\end{align*}
which implies that 
\begin{align*}
	\deg s_m(q)\bracq m 
	& = \deg s_m(q) + (m-1) < \deg u_m (q) + (m-1) ,
	\\
	\deg t_m(q)\bracq n 
	& = \deg u_m(q) + (n-1) ,
	\\
	\deg u_m(q)\bracq k
	& = \deg u_m(q) + (k-1) ,
	\\
	\deg v_m(q) \bracq n \bracq k
	& = \deg u_m(q) + (n-1) + (k-1) .
\end{align*}
In the case of $m<n+k-1$, we have
\begin{align*}
	\deg \bracq{m+n+k}
	& = 
	m + n + k -1 
	\\
	& < 
	\deg u_m(q) + (n-1) + (k-1) 
	= \deg v_m(q) \bracq n \bracq k
	\\
	& = 
	\deg \paren{ s_m (q) \bracq m + t_m (q) \bracq n + u_m (q) \bracq k + v_m (q) \bracq n \bracq k } ,
\end{align*}
a contradiction.
\endproof


\section{Zero identity}

For every double sequence of polynomial 
$\mathcal R = \set{r'_{m,n} (q)}_{m,n=1}^\infty$, 
and sequences of polynomials 
$\mathcal S = \set{s'_m (q)}_{m=1}^\infty$, 
$\mathcal T = \set{t'_m (q)}_{m=1}^\infty$, 
$\mathcal U = \set{u'_m (q)}_{m=1}^\infty$, 
$\mathcal V = \set{v'_m (q)}_{m=1}^\infty$, 
$\mathcal W = \set{w'_m (q)}_{m=1}^\infty$, 
define an addition rule for three $q$-integers by  
\begin{align*}
	\oplus_{\mathcal R, \, \mathcal U \mathcal S, \, \mathcal V \mathcal T, \mathcal W}
	= r'_{n,k} (q) \bracq m + u'_m (q) s'_k (q) \bracq n + v'_m (q) t'_n (q) \bracq k + w'_m (q) \bracq n \bracq k .
\end{align*}
This is called the first kind of quadratic zero identity for three $q$-integers, if 
\begin{align*}
	\oplus_{\mathcal R, \, \mathcal U \mathcal S, \, \mathcal V \mathcal T, \mathcal W}
	\paren{\bracq m, \bracq n, \bracq k} = 0 
\end{align*}
for all positive integers $m$, $n$ and $k$. 

\begin{theorem}
For every double sequence of polynomial 
$\mathcal R = \set{r'_{m,n} (q)}_{m,n=1}^\infty$, 
three sequences of polynomials 
$\mathcal S = \set{s'_m (q)}_{m=1}^\infty$, 
$\mathcal T = \set{t'_m (q)}_{m=1}^\infty$, 
$\mathcal W = \set{w'_m (q)}_{m=1}^\infty$, and
two sequences of polynomials 
$\mathcal U = \set{u'_m (q)}_{m=1}^\infty$ and
$\mathcal V = \set{v'_m (q)}_{m=1}^\infty$ 
with the intitial costants $u'_1(q)\equiv u$ and $v'_1(q) \equiv v$, respectively, the $\oplus_{\mathcal R, \, \mathcal U \mathcal S, \, \mathcal V \mathcal T, \mathcal W}$ is the first kind of quadratic zero identity for three $q$-integers if and only if
\begin{align*}
	r'_{n,k} (q)
	& = r_{1,k} (q) \bracq n + r_{n,1} (q) \bracq k - r(q) \bracq n \bracq k ,
	\\
	s'_k (q) 
	& = -\frac{1}{u} \set{ r_{1,k} (q) - r(q) \bracq k - u s(q) \bracq k}  ,
	\\
	t'_n (q) 
	& = -\frac{1}{v} \set{ r_{n,1} (q) - r(q) \bracq n - v t(q) \bracq n} ,
	\\
	u'_m (q) 
	& = u \bracq m ,
	\\
	v'_m (q) 
	& = v \bracq m ,
	\\
	w'_m (q) 
	& = - \set{ r(q) + u s(q) + v t(q) } \bracq m ,
\end{align*}
for some initial polynomials $s(q)$, $t(q)$, and two sequences of polynomials $\set{r_{1,k}(q)}_{k=1}^\infty$, $\set{r_{n,1}(q)}_{n=1}^\infty$ with a common intitial polynomial $r(q)=r_{1,1}(q)$.
\end{theorem}

\proof
Let $\oplus_{\mathcal R, \, \mathcal U \mathcal S, \, \mathcal V \mathcal T, \mathcal W}$ be the first kind of quadratic zero identity for three $q$-integers. Thus
\begin{align}
	r'_{n,k} (q) \bracq m + u'_m (q) s'_k (q) \bracq n + v'_m (q) t'_n (q) \bracq k + w'_m (q) \bracq n \bracq k = 0
\label{eq:zero:1}
\end{align}
for all positive integers $m,n,k$. 
Choose $r(q) = r'_{1,1}(q)$, $s(q) = s'_1(q)$ and $t(q) = t'_1(q)$ from the initial polynomails of the sequences $\mathcal R$, $\mathcal S$, $\mathcal T$, respectively.
For $m=1$, set $w (q)=w'_1 (q)$, we have that for all positive integers $n$ and $k$,
\begin{align}
	r'_{n,k} (q) \bracq 1 + u s'_k(q) \bracq n + v t'_n (q) \bracq k + w(q) \bracq n \bracq k =0.
\label{eq:zero:2}
\end{align}
For $m=1$ and $n=1$, choose $r_{1,k}(q) = r'_{1,k}(q)$, we obtain
\begin{align*}
	r_{1,k} (q) \bracq 1 + u s'_k(q) \bracq 1 + v t (q) \bracq k + w(q) \bracq 1 \bracq k =0,
\end{align*}
so that for all $k\ge1$,
\begin{align}
	s'_k(q) = \frac{-1}{u} \set{ r_{1,k} (q) + v t(q) \bracq k + w(q) \bracq k } 
\label{eq:zero:3}
\end{align}
For $m=1$ and $k=1$, choose $r_{n,1}(q) = r'_{n,1}(q)$, we also obtain
\begin{align*}
	r_{n,1} (q) \bracq 1 + u s(q) \bracq n + v t'_n (q) \bracq 1 + w(q) \bracq n \bracq 1 = 0,
\end{align*}
and so that for all $n\ge1$,
\begin{align}
	t'_n (q) = \frac{-1}{v} \set{ r_{n,1} (q) + u s(q) \bracq n + w(q) \bracq n } 
\label{eq:zero:4}
\end{align}
Replacing (\ref{eq:zero:3}) and (\ref{eq:zero:4}) in (\ref{eq:zero:2}), we can see that for all positive integers $n$ and $k$,
\begin{align}
	r'_{n,k} (q) 
	& = r_{1,k} (q) \bracq n + r_{n,1} (q) \bracq k
	+ \set{ u s(q) + v t (q) + w(q) } \bracq n \bracq k 
\label{eq:zero:5}
\end{align}
For $m=1$, $n=1$ and $k=1$, the quadratic zero identity (\ref{eq:zero:1}) can be reduced to 
\begin{align*}
	r(q) \bracq 1 + u s(q) \bracq 1 + v t(q) \bracq 1 + w(q) \bracq 1 \bracq 1 = 0,
\end{align*}
this implies that 
\begin{align}
	w(q) = - \set{r(q) + u s(q) + v t(q)} .
\label{eq:zero:6}
\end{align}
By replacing (\ref{eq:zero:6}) in (\ref{eq:zero:5}), (\ref{eq:zero:3}) and (\ref{eq:zero:4}), respectively, we obtain
\begin{align*}
	r'_{n,k} (q)
	& = r_{1,k} (q) \bracq n + r_{n,1} (q) \bracq k - r(q) \bracq n \bracq k ,
	\\
	s'_k (q) 
	& = -\frac{1}{u} \set{ r_{1,k} (q) - r(q) \bracq k - u s(q) \bracq k}  ,
	\\
	t'_n (q) 
	& = -\frac{1}{v} \set{ r_{n,1} (q) - r(q) \bracq n - v t(q) \bracq n} .
\end{align*}
From these relations, the quadratic zero identity (\ref{eq:zero:1}) can be expanded to 
\begin{align*}
	0= \;
	& r_{1,k} (q) \bracq n \bracq m 
	+ r_{n,1} (q) \bracq k \bracq m 
	- r(q) \bracq n \bracq k \bracq m 
	\\
	& - \frac{u'_m(q)}{u} r_{1,k} (q) \bracq n 
	+ \frac{u'_m(q)}{u} r(q) \bracq k \bracq n
	+ u'_m(q) s(q) \bracq k \bracq n
	\\
	& - \frac{v'_m(q)}{v} r_{n,1} (q) \bracq k 
	+ \frac{v'_m(q)}{v} r(q) \bracq n \bracq k
	+ v'_m(q) t(q) \bracq n \bracq k
\end{align*}
Since $r_{1,k}(q)$ and $r_{n,1}(q)$ are arbitrary polynomials, we can conclude that 
\begin{align*}
	0 = r_{1,k} (q) \bracq n \bracq m 
	- \frac{u'_m(q)}{u} r_{1,k} (q) \bracq n ,
	\\
	0 = r_{n,1} (q) \bracq k \bracq m 
	 - \frac{v'_m(q)}{v} r_{n,1} (q) \bracq k ,
\end{align*}
which implies that $u'_m(q) = u \bracq m$, and $v'_m(q) = v \bracq m$ for all positive integers $m$.
Finally, for $n=1$ and $k=1$, the quadratic zero identity (\ref{eq:zero:1}) can be reduced to 
\begin{align*}
	r (q) \bracq m + u'_m (q) s(q) \bracq 1 + v'_m(q) t (q) \bracq 1 + w'_m (q) \bracq 1 \bracq 1 = 0,
\end{align*}
and hence
\begin{align*}
	w'_m (q) 
	= - \set{r(q) + u s(q) + v t(q) } \bracq m .
\end{align*}

\endproof


\section{Polynomial functional equations}

In this section, we shall compute the solution of functional equations for a sequence of polynomials $\set{f_n(q)}_{n=1}^\infty$ associated to the first kind of quadratic addition rules:
\begin{align}
	f_{m+n+k} = f_m (q) + q^m f_n (q) + q^m f_k (q) + q^m (q-1) f_n (q) f_k (q) 
	\label{eq:quafe}
\end{align}
for all positive integers $m$, $n$ and $k$.
This functional equation always has the trivial solutions $f_n(q)=\bracq n$ for all $n\ge1$, and $f_n(q)\equiv0$ for all positive integers $n$.

\begin{theorem}
The polynomial functional equation (\ref{eq:quafe}) has only trivial solutions.
\end{theorem}

\proof
Suppose that the solution of (\ref{eq:quafe}) are given by $f_n (q) = h(q) \bracq n$ for some polynomial $h(q)$. Then 
\begin{align*}
	f_{n+2} (q) 
	&= f_{n+1+1} (q) 
	\\
	&= f_n (q) + q^n f_1 (q) + q^n f_1 (q) + q^n (q-1) f_1 (q) f_1(q)
	\\
	&= h(q) \bracq n + q^n h(q) \bracq 1 + q^n h(q) \bracq 1  + q^n (q-1) h(q) \bracq 1 h(q) \bracq 1 
	\\
	&=  h(q) \set{\bracq n + q^n \bracq 1 + q^n \bracq 1 + q^n (q-1) \bracq 1 \bracq 1 } 
	\\
	& \quad + q^n (q-1) \paren{h^2(q)-h(q)} \bracq 1 \bracq 1
	\\
	&= h(q) \bracq{n + 2} + q^n (q-1) \paren{h^2(q)-h(q)} \bracq 1 \bracq 1 
	\\
	&= f_{n+2}(q) + q^n (q-1) \paren{h^2(q)-h(q)} \bracq 1 \bracq 1 .
\end{align*}
So, that $0 = q^n (q-1) \paren{h^2(q)-h(q)} \bracq 1 \bracq 1$, 
and that $0 = h^2(q)-h(q)$.
\\
Hence, either $h(q)\equiv0$ or $h(q)\equiv1$.
This completes the proof.
\endproof


\section{The second kind of quadratic addition rules}

The functional equation for a sequence of polynomials $\set{f_n(q)}_{n=1}^\infty$ is given by 
\begin{align}
	f_{m+n} (q)
	&= f_m (q) + q^m f_n (q)
	\tag{linear functional eq.}
	\label{fe1}
	\\
	f_{m+n} (q)
	&= f_m (q) + f_n (q) + \paren{q-1} f_m (q) f_n (q)
	\tag{quadratic functional eq. $\mathrm{I}$}
	\label{fe2}
	\\
	f_{m+n} (q)
	&= q^n f_m (q) + q^m f_n (q) +(1-q) f_m (q) f_n (q)
	\tag{quadratic functional eq. $\mathrm{II}$}
	\label{fe3}
\end{align}


\bibliographystyle{chicago}
\bibliography{qFunctionalEquation}

\begin{thebibliography}{}

\bibitem[\protect\citeauthoryear{Kac and Cheung}{Kac and
  Cheung}{2001}]{kac2001quantum}
Kac, V. and P.~Cheung (2001).
\newblock {\em Quantum calculus}.
\newblock Springer Science \& Business Media.

\bibitem[\protect\citeauthoryear{Kassel}{Kassel}{2012}]{kassel2012quantum}
Kassel, C. (2012).
\newblock {\em Quantum groups}, Volume 155.
\newblock Springer Science \& Business Media.

\bibitem[\protect\citeauthoryear{Kontorovich and Nathanson}{Kontorovich and
  Nathanson}{2006}]{kontorovich2006quadratic}
Kontorovich, A.~V. and M.~B. Nathanson (2006).
\newblock Quadratic addition rules for quantum integers.
\newblock {\em Journal of Number Theory\/}~{\em 117\/}(1), 1--13.

\bibitem[\protect\citeauthoryear{Nathanson}{Nathanson}{2003}]{nathanson2003functional}
Nathanson, M.~B. (2003).
\newblock A functional equation arising from multiplication of quantum
  integers.
\newblock {\em Journal of Number Theory\/}~{\em 103\/}(2), 214--233.

\bibitem[\protect\citeauthoryear{Nathanson}{Nathanson}{2004}]{nathanson2004formal}
Nathanson, M.~B. (2004).
\newblock Formal power series arising from multiplication of quantum integers.
\newblock {\em DIMACs Series in Discrete Mathematics and Theoretical Computer
  Science\/}~{\em 64}, 145--168.

\bibitem[\protect\citeauthoryear{Nathanson}{Nathanson}{2007}]{nathanson2007linear}
Nathanson, M.~B. (2007).
\newblock Linear quantum addition rules.
\newblock {\em INTEGERS: ELECTRONIC JOURNAL OF COMBINATORIAL NUMBER
  THEORY\/}~{\em 7\/}(2), A27.

\bibitem[\protect\citeauthoryear{Nguyen}{Nguyen}{2010a}]{nguyen2010solutions}
Nguyen, L. (2010a).
\newblock On the solutions of a functional equation arising from multiplication
  of quantum integers.
\newblock {\em Journal of Number Theory\/}~{\em 130\/}(6), 1292--1347.

\bibitem[\protect\citeauthoryear{Nguyen}{Nguyen}{2010b}]{nguyen2010support}
Nguyen, L. (2010b).
\newblock On the support base of a functional equation arising from
  multiplication of quantum integers.
\newblock {\em Journal of Number Theory\/}~{\em 130\/}(6), 1348--1373.

\bibitem[\protect\citeauthoryear{Nguyen}{Nguyen}{2012}]{nguyen2012polynomial}
Nguyen, L. (2012).
\newblock On the polynomial and maximal solutions to a functional equation
  arising from multiplication of quantum integers.
\newblock {\em Notes on Number Theory and Discrete Mathematics\/}~{\em
  18\/}(4), 11--39.

\bibitem[\protect\citeauthoryear{Nguyen}{Nguyen}{2013}]{nguyen2013classification}
Nguyen, L. (2013).
\newblock On the classification of solutions of a functional equation arising
  from multiplication of quantum integers.
\newblock {\em Uniform Distribution Theory\/}~{\em 8\/}(2), 49--120.

\end{thebibliography}

%
%

\end{document}